\documentclass[a4paper]{article}
\usepackage{latexsym}
\usepackage{amsfonts,latexsym}
\usepackage{amsmath,amssymb}
\usepackage{amsthm}

\def\mod{{\rm mod\,}}

\newtheorem{Theorem}{Theorem}[section]
\newtheorem{Proposition}[Theorem]{Proposition}
\newtheorem{Lemma}[Theorem]{Lemma}

\newtheorem{Conjecture}[Theorem]{Conjecture}
\newcommand{\Proof}{{\bf Proof}\quad}
\newcommand{\C}{{\bf C}}
\newcommand{\nn}{{\mathcal N}}
\newcommand{\vv}{{\mathcal V}}
\newcommand{\SL}{\mathrm{SL}}
\newcommand{\doubleslash}{/\kern-1.5pt/}

\title{The degrees of a system of parameters of the ring
of invariants of a binary form}
\author{Andries E. Brouwer, Jan Draisma \& Mihaela Popoviciu}
\date{2009-08-18,~ 2014-04-18}

\begin{document}
\maketitle

\begin{abstract}
We consider the degrees of the elements of a homogeneous
system of parameters for the ring of invariants of a binary form,
give a divisibility condition, and a complete classification for
forms of degree at most 8.
\end{abstract}

\section{The degrees of a system of parameters}
Let $R$ be a graded $\C$-algebra.
A {\em homogeneous system of parameters} (hsop) of $R$ is an algebraically
independent set $S$ of homogeneous elements of $R$ such that $R$
is module-finite over the subalgebra generated by $S$.
By the Noether normalization lemma, a hsop always exists.
The size $|S|$ of $S$ equals the Krull dimension of $R$.

In this note we consider the special case where $R$ is
the ring $I$ of invariants of binary forms of degree $n$
under the action of $\SL(2,\C)$.
This ring is Cohen-Macaulay, that is, $I$ is free over the
subring generated by any hsop $S$. Its Krull dimension is $n-2$.

One cannot expect to classify all hsops of $I$. Indeed, any
generic subset with the right degrees will be a hsop
(cf.~Dixmier's criterion below). But one can expect to classify
the sets of degrees of hsops.
In this note we give a divisibility restriction on the set of degrees
for the elements of a hsop, and conjecture that when all degrees are
large this restriction also suffices for the existence of a hsop
with these given degrees. For small degrees there are further
restrictions. We give a complete classification for $n \le 8$.

\section{Hilbert's criterion}

Hilbert's criterion gives a characterization of
homogeneous systems of parameters as sets that define the nullcone.

Denote by $V_n$ the set of binary forms of degree $n$. The {\em
nullcone} of $V_n$, denoted $\nn(V_n)$, is the set of binary forms
of degree $n$ on which all invariants vanish. By the Hilbert-Mumford numerical
criterion (see \cite{Hi2} and \cite[Chapter 2]{MFK}) this is precisely
the set of binary forms of degree $n$ with a root of multiplicity
$>\frac{n}{2}$. Moreover, the binary forms with no root of multiplicity
$\geq \frac{n}{2}$ have closed $\SL(2,\C)$-orbits. The elements of
$\nn(V_n)$ are called {\em nullforms}. Another result from \cite{Hi2}
that we will use is the following.
\begin{Proposition} \label{hilbert} 
For $n \ge 3$, consider $i_1,\ldots ,i_{n-2}$ homogeneous invariants of $V_n$.
The following two conditions are equivalent:
\begin{itemize}
\item[(i)] $\nn(V_n)=\vv(i_1,\ldots ,i_{n-2})$,
\item[(ii)] $\{i_1,\ldots ,i_{n-2}\}$ is a hsop of the invariant ring of $V_n$.
\end{itemize}
\end{Proposition}

\section{A divisibility condition}
Assume $n \ge 3$.

\begin{Lemma} \label{lm:divisible}
Fix integers $j$, $t$ with $t > 0$.
If an invariant of degree $d$ is nonzero on a form $\sum a_i x^{n-i} y^i$
with the property that all nonzero $a_i$ have $i \equiv j$ (mod $t$),
then $d(n-2j)/2 \equiv 0$ (mod $t$).
\end{Lemma}
\Proof
For an invariant of degree $d$ with nonzero term $\prod a_i^{m_i}$ we have
$\sum m_i = d$ and $\sum i m_i = nd/2$.
If $i \equiv j$ (mod $t$) when $a_i \ne 0$, then
$nd/2 = \sum i m_i \equiv j \sum m_i = jd$ (mod $t$).
\qed

\medskip\noindent
For odd $n$ we recover the well-known fact that all degrees are even
(take $t=1$).

\begin{Lemma} Fix integers $j$, $t$ with $t > 1$ and $0 \le j \le n$.
Among the degrees $d$ of a hsop, at least
$\lfloor (n-j)/t \rfloor$ satisfy $d(n-2j)/2 \equiv 0$ (mod $t$).
\end{Lemma}
\Proof
Subtracting a multiple of $t$ from $j$ results in a stronger statement, so
it suffices to prove the lemma for $0 \le j < t$.  There are $1 + \lfloor
(n-j)/t \rfloor=:1+N$ coefficients $a_i$ with $i \equiv j$ (mod $t$),
so the subpace $U$ of $V_n$ defined by $a_i = 0$ for $i \not\equiv j$
(mod $t$) has dimension $1+N$. If $N=0$ there is nothing to prove, so
we assume that $N>0$. We claim that a general form $f \in U$ has only
zeroes of multiplicity strictly less than $n/2$. Indeed, write
\[ f=a_j x^{n-j} y^j + a_{j+t} x^{n-j-t} y^{j+t} + \ldots + 
a_{j+mt} x^{n-j-mt} y^{j+mt} \]
where $j+(m+1)t>n$ and $m>0$. So $f$ has a factor $y$ of multiplicity $j$
and a factor $x$ of multiplicity $n-j-mt$. If $j$ were at least $n/2$,
then $j+mt\geq j+t>2j\geq n$, a contradiction. If $n-j-mt$ were at least
$n/2$, then $j+mt \leq n/2$ and hence $t \leq n/2$ and hence $j+(m+1)t
\leq n$, a contradiction. The remaining roots of $f$ are roots of
\[ a_j x^{mt} + a_{j+t} x^{(m-1)t}y^t + \ldots + a_{j+mt}y^{mt}, \]
which is a general binary form of degree $m$ in $x^t,y^t$
and hence has $mt$ distinct roots. 

Let $\pi:V_n \to V_n\doubleslash\SL(2,\C)$ be the quotient map; so the
right-hand side is the spectrum of the invariant ring $I$.
Set $X:=\overline{\pi(U)}$. We claim that $X$ has dimension $N$.
It certainly cannot have dimension larger than $N$, since acting with the
one-dimensional torus of diagonal matrices on an element of $U$ gives
another element of $U$. To show that $\dim X=N$ we need to show that
for general $f \in U$ the fibre $\pi^{-1}(\pi(f))$ intersects $U$ in a
one-dimensional variety. By the above and the Hilbert-Mumford criterion,
the $\SL(2,\C)$-orbit of $f$ is closed. Moreover, its stabiliser is
zero-dimensional. So by properties of the quotient
map we have $\pi^{-1}(\pi(f))=\SL(2,\C) \cdot f$. Hence it suffices
that the intersection of this orbit with $U$ is one-dimensional. For
this a Lie algebra argument suffices, in which we may ignore the
Lie algebra of the torus: if $(b x \frac{\partial}{\partial y} + c y
\frac{\partial}{\partial x})f$ lies in $U$, then we find that $b=c=0$
if $t>2$ (so that the contribution of one term from $f$ cannot cancel the
contribution from the next term); and $b=0$ if $j>0$ (look at the first term),
and then also $c=0$; and $c=0$ if $j+mt<n$ (look at the last term),
and then also $b=0$. Hence the only case that remains is $t=2,j=0,$
and $n \geq 4$ even. Then the equations $c a_0 n + b a_2 2=0$ and $c a_2
(n-2) + b a_4 4=0$ are independent and force $b=c=0$. 

This concludes the proof that $\dim X=N$. Intersecting $X$ with the
hypersurfaces corresponding to elements of an hsop reduces $X$ to
the single point in $X$ representing the null-cone. In
the process, $\dim X$ drops by $N$. But the only invariants that
contribute to this dimension drop, i.e., the only invariants that do
not vanish identically on $X$ (hence on $U$) are those considered in
Lemma~\ref{lm:divisible}. Hence there must be at least $N$ of these
among the hsop.
\qed

\begin{Lemma}\label{lemma3}
Let $t$ be an integer with $t > 1$.

(i) If $n$ is odd, and $j$ is minimal such that $0 \le j \le n$ and
$(n-2j,t) = 1$, then among the degrees of any hsop at least
$\lfloor (n-j)/t \rfloor$ are divisible by $2t$.

(ii) If $n$ is even, and $j$ is minimal with $0 \le j \le \frac{1}{2}n$ and
$(\frac{1}{2}n-j,t) = 1$, then among the degrees of any hsop at least
$\lfloor (n-j)/t \rfloor$ are divisible by $t$. \qed
\end{Lemma}
%
% This is trivial. The nontrivial part is that this captures the
% full strength of the previous lemma.

\begin{Theorem}\label{thm1}
Let $t$ be an integer with $t > 1$.

(i) If $n$ is odd, then among the degrees of any hsop at least
$\lfloor (n-1)/t \rfloor$ are divisible by $2t$ (and all degrees are even).

(ii) If $n$ is even, then among the degrees of any hsop at least
$\lfloor (n-1)/t \rfloor$ are divisible by $t$, and if $n \equiv 2$ $(\mod 4)$
then at least $n/2$ by $2$.
\end{Theorem}
\Proof
(i) By part (i) of Lemma \ref{lemma3} we find a lower bound
$\lfloor (n-j)/t \rfloor$ for a $j$ as described there.
If that is smaller than $\lfloor (n-1)/t \rfloor$, then there is
some multple $at$ of $t$ with $n-j+1 \le at \le n-1$.
Put $n=at+b$, where $1 \le b \le j-1$. By definition of $j$ we have
$(b-2i,t) > 1$ for $i=0,1,...,j-1$. If $b$ is odd, say $b=2i+1$,
we find a contradiction. If $b$ is even, say $b=2i+2$, then $t$ is even
and $n$ is even, contradiction.

(ii) By part (ii) of Lemma \ref{lemma3} we find a lower bound
$\lfloor (n-j)/t \rfloor$ for a $j$ as described there. For $t = 2$
our claim follows. Now let $t > 2$.
If $\lfloor (n-j)/t \rfloor$ is smaller than $\lfloor (n-1)/t \rfloor$,
then there is some multple $at$ of $t$ with $n-j+1 \le at \le n-1$.
Put $n=at+b$, where $1 \le b \le j-1$. By definition of $j$ we have
$(b-2i,2t) > 2$ for $i=0,1,...,j-1$, impossible.
\qed

\medskip
For example, it is known that there exist homogeneous systems
of parameters with degree sequences
4 $(n=3)$;
2, 3 $(n=4)$;
4, 8, 12 $(n=5)$;
2, 4, 6, 10 $(n=6)$;
4, 8, 12, 12, 20 and 4, 8, 8, 12, 30 $(n=7)$ \cite{Di0};
2, 3, 4, 5, 6, 7 $(n=8)$ \cite{Sh};
4, 8, 10, 12, 12, 14, 16 and 4, 4, 10, 12, 14, 16, 24 and
4, 4, 8, 12, 14, 16, 30 and 4, 4, 8, 10, 12, 16, 42 and
4, 4, 8, 10, 12, 14, 48 $(n=9)$ \cite{nonic};
2, 4, 6, 6, 8, 9, 10, 14 $(n=10)$ \cite{decimic}.

\begin{Conjecture}
Any sequence $d_1,...,d_{n-2}$ of sufficiently large integers
satisfying the divisibility conditions of Theorem \ref{thm1}
is the sequence of degrees of a hsop.
\end{Conjecture}

This can be compared to the conjecture

\medskip
\begin{Conjecture} {\rm (Dixmier\cite{Di1})}

(i) If $n$ is odd, $n \ge 15$, then $4,6,8,...,2n-2$
is the sequence of degrees of a hsop.

(ii) If $n \equiv 2$ $(\mod 4)$, $n \ge 18$, then $2,4,5,6,6,7,8,9,...,n-1$
is the sequence of degrees of a hsop.

(iii) If $n \equiv 0$ $(\mod 4)$, then $2,3,4,...,n-1$
is the sequence of degrees of a hsop.
\end{Conjecture}

\section{Poincar\'e series}
If there exists a hsop with degrees $d_1,\ldots,d_{n-2}$,
then the Poincar\'e series can be written as a quotient
$P(t) = a(t) / \prod (t^{d_i}-1)$ for some polynomial $a(t)$
with nonnegative coefficients. If one does not have a hsop, but
only a sequence of degrees, the conditions of Theorem \ref{thm1} above
are strong enough to guarantee that $P(t)$ can be written in this way,
but without the condition that the numerator has nonnegative coefficients.

\begin{Proposition}
Let $d_1,\ldots,d_{n-2}$ be a sequence of positive integers
satisfying the conditions of Theorem \ref{thm1}.
Then $P(t) \prod (t^{d_i}-1)$ is a polynomial.
\end{Proposition}
\Proof
Dixmier \cite{Di1} proves that $P(t) B(t)$ is a polynomial,
where $B(t)$ is defined by
\[
B(t) = \left\{ \begin{array}{ll}
\prod_{i=2}^{n-1} (1-t^{2i}) & \mbox{if $n$ is odd} \\[2pt]
\prod_{i=2}^{n-1} (1-t^i).(1+t) & \mbox{if $n \equiv 2$ (mod 4)} \\[2pt]
\prod_{i=2}^{n-3} (1-t^i).(1+t)(1-t^{(n-2)/2})(1-t^{n-1}) &
\mbox{if $n \equiv 0$ (mod 4)} \\
\end{array} \right.
\]
Consider a primitive $t$-th root of unity $\zeta$.
We have to show that if $B(t)$ has root $\zeta$ with multiplicity $m$,
then at least $m$ of the $d_i$ are divisible by $t$,
but this follows immediately from Theorem \ref{thm1}.
Note that in case $n \equiv 0$ (mod 4) the factor $(1+t)(1-t^{(n-2)/2})$
divides $(1-t^{n-2})$.
\qed

\medskip\noindent
We see that if $n \equiv 0$ (mod 4), $n > 4$, then $P(t)$ can be written
with a smaller denominator than corresponds to the degrees of a hsop.

\begin{table}[ht]
\medskip\noindent\begin{center}{\small
\begin{tabular}{@{}r|rrrrrrrrrrr@{~}r@{~}r@{~}r@{~}r@{}}
$h^n_m$ & 1 & 2 & 3 & 4 & 5 & 6 & 7 & 8 & 9 & 10 & 11 & 12 & 13 & 14 & 15 \\
\hline
1 & . & . & . & . & . & . & . & . & . & . & . & . & . & . & . \\
2 & . & 1 & . & 1 & . & 1 & . & 1 & . & 1 & . & 1 & . & 1 & .\\
3 & . & . & . & 1 & . & . & . & 1 & . & . & . & 1 & . & . & .\\
4 & . & 1 & 1 & 1 & 1 & 2 & 1 & 2 & 2 & 2 & 2 & 3 & 2 & 3 & 3\\
5 & . & . & . & 1 & . & . & . & 2 & . & . & . & 3 & . & . & .\\
6 & . & 1 & . & 2 & . & 3 & . & 4 & . & 6 & . & 8 & . & 10 & 1\\
7 & . & . & . & 1 & . & . & . & 4 & . & . & . & 10 & . & 4 & .\\
8 & . & 1 & 1 & 2 & 2 & 4 & 4 & 7 & 8 & 12 & 13 & 20 & 22 & 31 & 36\\
9 & . & . & . & 2 & . & . & . & 8 & . & 5 & . & 28 & . & 27 & .\\
10 & . & 1 & . & 2 & . & 6 & . & 12 & 5 & 24 & 13 & 52 & 33 & 97 & 80\\
11 & . & . & . & 2 & . & . & . & 13 & . & 13 & . & 73 & . & 110 & .\\
12 & . & 1 & 1 & 3 & 3 & 8 & 10 & 20 & 28 & 52 & 73 & 127 & 181 & 291 & 418\\
13 & . & . & . & 2 & . & . & . & 22 & . & 33 & . & 181 & . & 375 & .\\
14 & . & 1 & . & 3 & . & 10 & 4 & 31 & 27 & 97 & 110 & 291 & 375 & 802 & 1111\\
15 & . & . & . & 3 & . & 1 & . & 36 & . & 80 & . & 418 & . & 1111 & .\\
16 & . & 1 & 1 & 3 & 4 & 13 & 18 & 47 & 84 & 177 & 320 & 639 & 1120 & 2077 & 3581 \\
17 & . & . & . & 3 & . & 1 & . & 54 & . & 160 & . & 902 & . & 2930 & . \\
18 & . & 1 & . & 4 & 1 & 16 & 13 & 71 & 99 & 319 & 529 & 1330 & 2342 & 5034 & 8899 \\
\end{tabular}}\end{center}
\caption{Values of $h^n_m = \dim_\C I_m$ with $I$ the ring of invariants
of a binary form of degree $n$. Here . denotes 0. One has $h^n_m = h^m_n$
and $P(t) = \sum_m h^n_m t^m$.
\label{tab1}}
\end{table}

\medskip
We shall need the first few coefficients of $P(t)$.
Messy details arise for small $n$ because
there are too few invariants of certain small degrees.
Let $I$ be the ring of invariants of a binary form of degree (order) $n$,
let $I_m$ be the graded part of $I$ of degree $m$,
and put $h_m = h^n_m = \dim_\C I_m$,
so that $P(t) = \sum_m h_m t^m$.

The coefficients $h^n_m$ can be computed by the Cayley-Sylvester formula:
The dimension of the space of covariants of degree $m$ and order $a$
is zero when $mn-a$ is odd, and equals $N(n,m,t)-N(n,m,t-1)$
if $nm-a = 2t$, where $N(n,m,t)$ is the number of ways $t$ can be
written as sum of $m$ integers in the range $0..n$, that is,
the number of Ferrers diagrams of size $t$ that fit into a $m \times n$
rectangle.

We have Hermite reciprocity $h^n_m = h^m_n$, as follows immediately
since reflection in the main diagonal shows $N(n,m,t) = N(m,n,t)$.
That means that Table \ref{tab1} is symmetric.

\medskip
Dixmier \cite{Di1} gives the cases in which $h_m = 0$.
Since his statement is not precisely accurate, we repeat his proof.

\begin{Proposition}\label{zeroh}
Let $m,n \geq 1$.
One has $h_m = h^n_m = 0$ precisely in the following cases:
\begin{itemize}
\item[(i)] if $mn$ is odd,
\item[(ii)] if $m=1$; if $n = 1$,
\item[(iii)] if $m=2$ and $n$ is odd; if $n=2$ and $m$ is odd,
\item[(iv)] if $m = 3$ and $n \equiv 2$ $(\mod~4)$;
if $n=3$ and $m \equiv 2$ $(\mod~4)$,
\item[(v)] if $m = 5$ and $n = 6,10,14$;
if $n = 5$ and $m = 6,10,14$,
\item[(vi)] if $m = 6$ and $n = 7,9,11,13$;
if $n = 6$ and $m = 7,9,11,13$,
\item[(vii)]
if $m = 7$ and $n = 10$;
if $n = 7$ and $m = 10$.
\end{itemize}
\end{Proposition}
\Proof
(i) If $n$ is odd, then all degrees are even.
(ii) For $n=1$ we have $P(t) = 1$.
(iii) For $n=2$ we have $P(t) = 1/(1-t^2)$.
(iv) For $n=3$ we have $P(t) = 1/(1-t^4)$.
Now let $m,n \ge 4$.
For $n = 4$ we have invariants of degrees 2, 3 and hence of all degrees
$m \ne 1$. That means that $h^n_4 \ne 0$.
For $n = 6$ we have invariants of degrees 2, 15 and hence of all degrees
$m \ge 14$. That means that $h^n_6 \ne 0$ for $n \ge 14$.
If $n$ is odd this shows the presence of invariants of degrees 4, 6
and hence of all even degrees $m > 2$, provided $n \ge 15$.
For $n = 5$ we have invariants of degrees 4, 18 and hence of all even degrees
$m \ge 16$. That means that $h^n_5 \ne 0$ for even $n \ge 16$.
If $n$ is even this shows the presence of invariants of degrees 2, 5
and hence of all degrees $m \ge 4$, provided $n \ge 16$.
It remains only to inspect the table for $4 \le m,n \le 14$.
\qed

\section{Dixmier's criterion}
Dividing out the ideal spanned by $p$ elements of a hsop
diminishes the dimension by precisely (and hence at least) $p$.
This means that the below gives a necessary and sufficient
condition for a sequence of degrees to be the degree sequence of a hsop.

\begin{Proposition} {\rm (Dixmier \cite{Di1})}
Let $G$ be a reductive group over $\C$, with a rational representation
in a vector space $R$ of finite dimension over $\C$. Let $\C[R]$ be the
algebra of complex polynomials on $R$, $\C[R]^G$ the subalgebra of
$G$-invariants, and $\C[R]^G_d$ the subset of homogeneous polynomials
of degree $d$ in $\C[R]^G$. Let $V$ be the affine variety such that
$\C[V] = \C[R]^G$. Let $r = \dim V$. Let $(d_1,\ldots,d_r)$
be a sequence of positive integers. Assume that for each subsequence
$(j_1,\ldots,j_p)$ of $(d_1,\ldots,d_r)$ the subset
of points of $V$ where all elements of all $\C[R]^G_j$ with
$j \in \{j_1,\ldots,j_p\}$ vanish has codimension not less than $p$
in $V$. Then $\C[R]^G$ has a system of parameters of degrees
$d_1,\ldots,d_r$. \qed
\end{Proposition}

This criterion is very convenient, it means that one can work with
degrees only, without worrying about individual elements of a hsop.

\section{Minimal degree sequences}
If $y_1,...,y_r$ is a hsop, then also $y_1^{e_1},...,y_r^{e_r}$
for any sequence of positive integers $e_1,...,e_r$, not all 1.
This means that if the degree sequence $d_1,...,d_r$ occurs,
also the sequence $d_1e_1,...,d_re_r$ occurs.
We would like to describe the minimal sequences, where such multiples
are discarded.

There are further reasons for non-minimality.

\begin{Lemma}\label{mainlemma}
If there exist hsops with degree sequences
$d_1,...,d_{r-1},d'$ and $d_1,...$, $d_{r-1}, d''$, then
there also exists a hsop with degree sequence $d_1,...,d_{r-1},d'+d''$.
\end{Lemma}
\Proof
We verify Dixmier's criterion. Consider a finite basis $f_1,...,f_s$
for the space of invariants of degree $d'$. Split the variety $V$
in the $s$ pieces defined by $f_i \ne 0$ $(1 \le i \le s)$
together with the single piece defined by $f_1 = ... = f_s = 0$.
Given $p$ elements of the sequence $d_1,...,d_{r-1},d'+d''$
we have to show that the codimension in $V$ obtained by requiring
all invariants of such degrees to vanish is at least $p$, that is,
that the dimension is at most $r-p$.
This is true by assumption if $d'+d''$ is not among these $p$ elements.
Otherwise, consider the $s+1$ pieces separately. We wish to show
that each has dimension at most $r-p$, then the same will hold for
their union. For the last piece, where all invariants of degree $d'$
vanish, this is true by assumption. But if some invariant of degree $d'$
does not vanish, and all invariants of degree $d'+d''$ vanish, then
all invariants of degree $d''$ vanish, and we are done.
\qed

\medskip
Note that taking multiples is a special case of (repeated application of)
this lemma, used with $d' = d''$.

Let us call a sequence {\em minimal} if it occurs (as the degree sequence
of the elements of a hsop), and its occurrence is not a consequence,
via the above lemma or via taking multiples,
of the occurrence of smaller sequences.
We might try to classify all minimal sequences, at least in small cases.

\medskip
Is it perhaps true that a hsop exists for any degree sequence
that satisfies the conditions of Theorem \ref{thm1} when there are
sufficiently many invariants? E.g. when the coefficients of
$P(t) \prod (1-t^{d_i})$ are nonnegative?

\medskip\noindent
{\bf Example}~
Some caution is required. For example, look at $n=6$. The conditions
of Theorem \ref{thm1} are: at least three factors 2, at least one factor
of each of 3, 4, 5. The sequence 6, 6, 6, 20 satisfies this restriction.
Moreover, $P(t) (1-t^6)^3 (1-t^{20}) =
1 + t^2 + 2t^4 + t^8 + 2t^{12} + t^{14} + t^{15} + t^{16} + t^{17} +
2t^{19} + t^{23} + 2t^{27} + t^{29} + t^{31}$
has only nonnegative coefficients.
But no hsop with these degrees exists: since $h_2 = 1$, $h_4 = 2$, $h_6 = 3$
it follows that there are invariants $i_2, i_4, i_6$ of degrees 2, 4, 6,
and we have $I_4 = \langle i_2^2, i_4 \rangle$ and
$I_6 = \langle i_2^3, i_2i_4, i_6 \rangle$. Requiring all invariants
of degree 6 to vanish is equivalent to the two conditions $i_2 = i_6 = 0$,
and a hsop cannot contain three elements of degree 6.

\medskip
Still, the above conditions almost suffice. And for $n < 6$ they
actually do suffice.

\subsection{$n=3$}

For $n=3$ we only have simple multiples of the minimal degree.

\begin{Proposition}
A positive integer $d$ is the degree of a hsop in case $n = 3$
if and only if it is divisible by $4$. \qed
\end{Proposition}

\noindent
If $i_4$ is an invariant of degree 4, then $\{i_4\}$ is a hsop.
% If $f\in V_3$ and $i_4=((f,f)_2,(f,f)_2)_2$, then $\{i_4\}$ is a hsop. 
% (Here $(g,h)_r$ denotes the $r$-th order transvectant of the forms
% $g$ and $h$, cf.~\cite{Olver}, p.88.)

\subsection{$n=4$}

For $n=4$ one has the sequence 2, 3, but for example also 5, 6.

\begin{Proposition}
A sequence $d_1,d_2$ of two positive integers is the sequence
of degrees of a hsop for the quartic if and only if neither of
them equals $1$, at least one is divisible by $2$, and at least one
is divisible by $3$.
\end{Proposition}
\Proof
Clearly the conditions are necessary.
In order to show that they suffice apply induction
and the known existence of a hsop with degrees 2, 3.
If $d_2 > 7$, then apply Lemma \ref{mainlemma} to
the two sequences $d_1,6$ and $d_1,d_2-6$ to conclude
the existence of a hsop with degrees $d_1,d_2$.
If $2 \le d_1,d_2 \le 7$ and one is divisible by 2,
the other by 3, then we have a multiple of the sequence 2, 3.
Otherwise, one equals 6 and the other is 5 or 7.
But 5, 6 is obtained from 2, 6 and 3, 6,
and 7, 6 is obtained from 2, 6 and 5, 6.
\qed

\medskip\noindent
If $i_2$ and $i_3$ are invariants of degrees 2 and 3,
then $\{i_2,i_3\}$ is a hsop.
% If $f\in V_4$ and $i_2=(f,f)_4$, $i_3=((f,f)_2,f)_4$,
% then $\{i_2,i_3\}$ is a hsop.  

\begin{Proposition}
There is precisely one minimal degree sequence of hsops in case $n = 4$,
namely $2$, $3$. \qed
\end{Proposition}

\subsection{$n=5$}

\begin{Proposition}
A sequence $d_1,d_2,d_3$ of three positive integers is the sequence
of degrees of a hsop for the quintic if and only if 
all $d_i$ are even, and distinct from $2$, $6$, $10$, $14$,
and no two are $4$, $4$ or $4$, $22$ and
at least two are divisible by $4$, at least one is divisible by $6$,
and at least one is divisible by $8$.
\end{Proposition}
\Proof
For $n=5$ the Poincar\'e series is
$P(t) = 1 + t^4 + 2t^8 + 3t^{12} + 4t^{16} + t^{18} + 5t^{20} +
t^{22} + 7t^{24} + 2t^{26} + 8t^{28} + 3t^{30} + ...$.
The stated conditions are necessary: the divisibility conditions
are seen from Theorem \ref{thm1}, and there are no invariants
of degrees 2, 6, 10, 14. Finally, we have $h_4 = 1$ and
$h_{18} = h_{22} = 1$, so that there are unique invariants $i_4$ and $i_{18}$
of degrees 4 and 18, respectively, and $I_{22} = \langle i_4i_{18} \rangle$,
so that all invariants of degree 22 will vanish as soon as $i_4$ vanishes.

The stated conditions suffice:
We use (and verify below) that there are hsops with degrees 4, 8, 12 and
with degrees 4, 8, 18.
If all $d_i$ are divisible by 4, and we do not have a multiple of 4, 8, 12,
then we have $4a$, $4b$, $24c$ where $a$ and $b$ have no factor 2 or 3,
and not both are 1.
It suffices to find 4, $4b$, 24. Since 4, 8, 24 exists, we can
decrease $b$ by 2, and it suffices to find 4, 12, 24, which exists.

So, some $d_i$, is not divisible by 4. We have one of the three cases
$24a,4b,2c$ and $8a,12b,2c$ and $8a,4b,6c$, where $c$ is odd.
In the middle case we have $c \ge 9$ and it suffices to make $8,12,2c$.
Since 8, 12, 4 exists, we can reduce $c$ by 2, and it suffices to make
8, 12, 18, which exists since 4, 8, 18 exists.

In the first case we have $c \ge 9$ and it suffices to make $24,4,2c$.
Since 12, 4, 8 exists, we can reduce $c$ by 4, and it suffices to make
24, 4, 18 and 24, 4, 30. The former is a multiple of 4, 8, 18 and the latter
follows from 24, 4, 18 and 24, 4, 12. Since 24, 4, 22 does not exist we still
have to consider $24a,4b,22$. Since 8, 12, 22 exists we can reduce $b$
by 2, and it suffices to make 24, 12, 22. But that is a multiple of 8, 12, 22.

Finally in the last case we have $c \ge 3$, and since 8, 4, 12 exists
we can reduce $c$ by 2. So it suffices to do 4, 8, 18, and that exists.
\qed

\begin{Proposition}
There are precisely two minimal degree sequences of hsops in case $n = 5$,
namely $4,8,12$ and $4,8,18$.
\end{Proposition}
\Proof
By the proof of the previous proposition, all we have to do is show
the existence of hsops with the indicated degree sequences.
It is well-known (see, e.g., Schur \cite{Schur}, p.86)
that the quintic has four invariants $i_4$, $i_8$, $i_{12}$, $i_{18}$
(with degrees as indicated by the index) that generate the ring of invariants,
and every invariant of degree divisible by 4 (in particular $i_{18}^2$)
is a polynomial in the first three.
Thus, when $i_4$, $i_8$, $i_{12}$ vanish, all invariants vanish,
and $\{i_4,i_8,i_{12}\}$ is a hsop.
Knowing this, it is easy to see that also $\{i_4,i_8,i_{18}\}$ is a hsop:
a simple Groebner computation shows that $i_{12}^3 \in (i_4,i_8,i_{18})$,
hence $\nn(V_5)=\vv (i_4,i_8,i_{18})$.
\qed

\subsection{$n=6$}

Similarly, we find for $n = 6$:

\begin{Proposition}
A sequence $d_1,d_2,d_3,d_4$ of four positive integers is the sequence
of degrees of a hsop for the sextic if and only if all $d_i$ are
distinct from $1$, $3$, $5$, $7$, $9$, $11$, $13$, and no two are
in $\{2,17\}$, and no three are in $\{2,4,8,14,17,19,23,29\}$,
and no three are in $\{2,6,17,21\}$,
and at least three are divisible by $2$,
at least one is divisible by $3$,
at least one by $4$, and at least one by $5$.
\end{Proposition}
\Proof
For $n=6$ the Poincar\'e series is
\begin{eqnarray*}
P(t) & = & 1 + t^2 + 2t^4 + 3t^6 + 4t^8 + 6t^{10} + 8t^{12} +
10t^{14} + t^{15} + 13t^{16} + t^{17} + \\
&& 16t^{18} + 2t^{19} + 20t^{20} + 3t^{21} +
24t^{22} + 4t^{23} + 29t^{24} + 6t^{25} + 34t^{26} + \\
&& 8t^{27} + 40t^{28} + 10t^{29} + 47t^{30} + \cdots .
\end{eqnarray*}
We have
\[
I_2 = \langle i_2 \rangle,~~
I_4 = \langle i_2^2, i_4 \rangle,~~
I_6 = \langle i_2^3, i_2i_4, i_6 \rangle,~~
I_8 = \langle i_2^4, i_2^2i_4, i_2i_6, i_4^2 \rangle,
\]
\[
I_{10} = \langle i_2^5, i_2^3i_4, i_2^2i_6, i_2i_4^2, i_4i_6, i_{10} \rangle,~~
I_{12} = \langle i_2^6, i_2^4i_4, i_2^3i_6, i_2^2i_4^2, i_2i_4i_6,
i_2i_{10}, i_4^3, i_6^2 \rangle,
\]
\[
I_{14} = \langle i_2^7, i_2^5i_4, i_2^4i_6, i_2^3i_4^2, i_2^2i_4i_6,
i_2^2i_{10}, i_2i_4^3, i_2i_6^2, i_4^2i_6, i_4i_{10} \rangle,~~
I_{15} = \langle i_{15} \rangle,
\]
and the invariants in degrees 17, 19, 23, 29 are $i_{15}$ times
the invariants in degrees 2, 4, 8, 14, respectively.
Let us denote by $[i_1,...,i_t]$ the condition that all invariants
of degrees $i_1, ..., i_t$ vanish. Then $[2] = [2,17]$ and hence
a hsop cannot have two element degrees among 2, 17.
Also $[4] = [2,4,8,14,17,19,23,29]$ and hence a hsop cannot have three
element degrees among 2, 4, 8, 14, 17, 19, 23, 29.
And $[6] = [2,6,17,21]$ is the condition $i_2 = i_6 = 0$ so that a hsop
cannot have three element degrees among 2, 6, 17, 21.
It follows that the stated conditions are necessary.

\medskip
The stated conditions suffice:
We use (and verify below) that there are hsops with each of the
degree sequences 2, 4, 6, 10 and 2, 4, 6, 15 and 2, 4, 10, 15.
Prove by induction that any 4-tuple of degrees that satisfies
the given conditions occurs as the degree sequence of a hsop.
Given $d_1,d_2,d_3,d_4$, if $d_i \geq 90$ then by induction we
already have the 4-tuples obtained by replacing $d_i$ by 60
and by $d_i-60$. It remains to check the finitely many cases
where all $d_i$ are less than 90. A small computer check settles this.
\qed

\begin{Proposition}
There are precisely three minimal degree sequences of hsops in case $n = 6$,
namely $2,4,6,10$ and $2,4,6,15$ and $2,4,10,15$.
\end{Proposition}
\Proof % MP
By the proof of the previous proposition, all we have to do is show
the existence of hsops with the indicated degree sequences.
It is well-known (see, e.g., Schur \cite{Schur}, p.90)
that the sextic has five invariants $i_2$, $i_4$, $i_6$, $i_{10}$, $i_{15}$
(with degrees as indicated by the index) that generate the ring of invariants,
where $i_{15}^2$ is a polynomial in the first four. This
implies that $\nn(V_6)=\vv(i_2,i_4,i_6,i_{10})$, so that
$\{i_2,i_4,i_6,i_{10}\}$ is a hsop.     
Now $\{i_2,i_4,i_6,i_{15}\}$ and $\{i_2,i_4,i_{10},i_{15}\}$ are also hsops:
we verified by computer that $i_{10}^3 \in (i_2,i_4,i_6,i_{15})$ and
$i_6^5\in (i_2,i_4,i_{10},i_{15})$, so that
$\nn(V_6)=\vv(i_2,i_4,i_6,i_{15})=\vv(i_2,i_4,i_{10},i_{15})$.
\qed

\subsection{$n=7$}

For $n = 7$ we have to consider the invariants a bit more closely
in order to decide which degree sequences are admissable for hsops.

Let $f$ be our septimic and let $\psi$ be the covariant $\psi = (f,f)_6$.
There are thirty basic invariants, of degrees 4, 8 (3$\times$),
12 (6$\times$), 14 (4$\times$), 16 (2$\times$), 18 (9$\times$),
20, 22 (2$\times$), 26, 30.
These can all be taken to be transvectants with a power of $\psi$
except for three basic invariants of degrees 12, 20 and 30
(that von Gall \cite{vG} calls $R$, $A$, $B$ and Dixmier \cite{Di0}
$q_{12}$, $p_{20}$, $p_{30}$).
This means that all invariants of degrees not of the form $12a+20b+30c$
vanish on the set defined by $\psi = 0$. But $\psi$ is a covariant of order 2,
i.e., $\psi = Ax^2 + Bxy + Cy^2$ for certain $A$, $B$, $C$.
It follows that no hsop degree sequence can have four elements in the set
$\{4,8,14,16,18,22,26,28,34,38,46,58\}$.

\begin{Proposition}
A sequence of five positive even integers is the sequence of degrees
of a hsop for the septimic if and only if all are distinct from 
$2$, $6$, $10$, no two equal $4$, no four are in
$\{4,8,14,16,18,22,26,28,34,38,46,58\}$
and at least three are divisible by $4$, at least two by $6$,
at least one by $8$, at least one by $10$ and at least one by $12$.
\end{Proposition}
\Proof
We already saw that these conditions are necessary.
For sufficiency, use induction.
The divisibility conditions concern moduli with l.c.m. 120,
and the restrictions concern numbers smaller than 60,
so if one of the degrees is not less than 180, we are done by induction.
A small computer program checks all degree sequences with degrees
at most 180, and finds that all can be reduced to the 23 sequences
given in the following proposition.
\qed

\begin{Proposition}
There are precisely $23$ minimal degree sequences of hsops in case $n = 7$,
namely

\[
\begin{array}{llll}
4,8,8,12,30    & 4,12,12,12,40  & 4,12,18,18,40 & 8,12,12,14,20 \\
4,8,12,12,20   & 4,12,12,14,40  & 4,14,14,24,60 & 8,12,14,14,60 \\
4,8,12,12,30   & 4,12,12,18,40  & 4,14,18,20,24 & 8,12,14,18,20 \\
4,8,12,14,30   & 4,12,14,14,120 & 4,14,18,32,60 & 12,12,14,14,40 \\
4,8,12,18,20   & 4,12,14,18,40  & 4,18,18,20,24 & 12,14,14,20,24 \\
4,8,12,18,30   & 4,12,14,20,24  & 4,18,18,32,60 \\
\end{array}
\]
\end{Proposition}
\Proof
We only have to show existence. Apply Dixmier's criterion.
Denote by $[d_1,...,d_p]$ the codimension in $V$ of the subset
of points of $V$ where all elements of all $\C[R]^G_{d_j}$
vanish $(1 \le j \le p)$. We have to show that for all $p$ and each of
these 23 sequences $(d_i)$ the inequality $[d_1,...,d_p] \ge p$ holds.

For $p=1$ that means that we need $[m] \ge 1$ for
$m = 4$, 8, 12, 14, 18, 20, 24, 30, 32, 40, 60, 120, and that is true,
for example by inspection of Table \ref{tab1}.

We can save some work by observing that Dixmier \cite{Di0} already showed
the existence of hsops with degree sequences 4, 8, 8, 12, 30 and
4, 8, 12, 12, 20. It follows that $[8] \ge 3$ and $[12] \ge 3$ and
$[24] \ge [8,12] \ge 4$ and $[20] \ge 2$ and $[60] \ge [12,20] \ge 4$
and $[4,30] \ge 2$ and $[8,30] \ge 4$.
Since there are several basic invariants of degree 14 or 18,
no two of which can have a common factor, it follows that
$[14] \ge 2$ and $[18] \ge 2$. This suffices to settle $p=2$.

For $p=3$ we must look at triples $[d,d',d'']$ without element
8 or 12 or multiple. First check that $[4,14] \ge 3$ and $[4,18] \ge 3$.
We'll do this below. Now all the rest needed for $p=3$ follows.

Below we shall show that $[12] \ge 4$.
For $p=4$ we must look at quadruples $[d,d',d'',d''']$ without
element 12 or 8, 30 or multiple.
The minimal of these are (omitting implied elements)
$[18,20]$ and $[18,32]$. However, $[18,32] \ge \min([18,12],[18,20])$
and $[18,20] \ge \min([18,20,8],[18,20,12])$.

Finally for $p=5$ we have to show that each of these 23 sets
determines the nullcone. But that follows immediately, since
it is known already that $[8,12,20] = [8,12,30] = 5$.

Altogether, our obligations are:
show that $[4,14] \ge 3$, $[4,18] \ge 3$, $[12] \ge 4$ and $[8,18,20] \ge 4$.

\medskip
Consider the part of $V$ defined by $\psi=0$.
Dixmier shows that if $\psi=q_{12}=p_{20}=0$ (for certain invariants
$q_{12}$ and $p_{20}$ of degrees 12 and 20, respectively), then
$f$ is a nullform. It follows that the subsets of $V$ defined by
$\psi=q_{12}=0$ or by $\psi=p_{20}=0$ have codimension at least 4 in $V$.

Now we have to do some actual computations.
With $f = ax^7+\binom{7}{1}bx^6y+\cdots+\binom{7}{1}gxy^6+hy^7$
(the two meanings of $f$, as form and as coefficient will not
cause confusion), we find
$\psi = (ag-6bf+15ce-10d^2)x^2 + (ah-5bg+9cf-5de)xy + (bh-6cg+15df-10e^2)y^2$.

Assume that the invariant of degree 4 vanishes, as it does in all cases
we still have to consider. Then $\psi$ has zero discriminant.
If $\psi \ne 0$, then w.l.o.g. $\psi \sim x^2$, and
$ah-5bg+9cf-5de = bh-6cg+15df-10e^2 = 0$, $ag-6bf+15ce-10d^2 \ne 0$.

Distinguish the four cases (i) $h \ne 0$, (ii) $h = 0$, $g \ne 0$,
(iii) $h = g = 0$, $f \ne 0$, (iv) $h = g = f = 0$, $e \ne 0$.
W.l.o.g. these become
(i) $h = 1$, $g = 0$, $a+9cf-5de = 0$, $b+15df-10e^2 = 0$,
(ii) $h = 0$, $g = 1$, $f = 0$, $b+de = 0$, $3c+5e^2 = 0$,
(iii) $h = g = 0$, $f = 1$, $e = 0$, $c = 0$, $d = 0$, $b \ne 0$,
(iv) $h = g = f = 0$, $e = 1$, $d = 0$, contradiction.

\medskip
Let us first show that $[12] \ge 4$. We may suppose $\psi \ne 0$.
One of the invariants of degree 12 is
$(\psi_1,\psi^5)_{10} \sim (\psi_1,x^{10})_{10} = fh-g^2$,
where $\psi_1 = (f,f)_2$.
If all invariants of degree 12 vanish, then in case (i) $f=0$,
and in case (ii) contradiction.
Look at case (iii). The only invariant of degree 12 that does not
vanish identically is $a^2b^2f^8$, and we find $a = 0$, a 1-dimensional set.
Finally, in case (i), if all invariants of degree 12 vanish,
but $ag-6bf+15ce-10d^2 \ne 0$, then the remaining conditions
define an ideal $(18e^3-cd,12de^2-c^2,2cd^2-3c^2e)$
in the three variables $c,d,e$ and the quotient is 1-dimensional.
This shows that $[12] \ge 4$.

\medskip
Let us show next that $[8,18] \ge 4$. We may suppose $\psi \ne 0$.
One of the invariants of degree 8 is
$(\psi_2,\psi^3)_6 \sim (\psi_2,x^6)_6 = dh-4eg+3f^2$
where $\psi_2 = (f,f)_4$.
This gives a contradiction in case (iii).
In case (ii) it gives $e=b=c=0$, leaving only variables $a,d$.
In case (i) it gives $d+3f^2=0$, leaving only variables $c,e,f$.

An invariant of degree 18 is
$((\psi_1,\psi_2)_1,\psi^7)_{14}  \sim ((\psi_1,\psi_2)_1,x^{14})_{14} =
-cfh^2+cg^2h+deh^2+2dfgh-3dg^3-4e^2gh+ef^2h+6efg^2-3f^3g$.
In case (ii) this says $d = 0$, leaving only variable $a$.
In case (i) this says $f(2ef+c) = 0$.
This gives us two subcases: (ia) with $f=0$ and variables $c,e$,
and (ib) with $c+2ef=0$ and variables $e,f$.

Another invariant of degree 8 is
$(\psi_3,\psi^2)_4 \sim (\psi_3,x^4)_4$, where $\psi_3 = (\psi_2,\psi_2)_4$,
which vanishes in case (ii) and says $c^2f+4cef^2+76e^2f^3+9e^4+144f^6 = 0$
in case (i). In case (ia) this means $e = 0$ leaving only variable $c$.
In case (ib) this means $(4f^3+e^2)^2 = 0$, leaving the dimension 1.
This proves $[8,18] \ge 4$.

\medskip
Let us show next that $[4,14] \ge 3$.
First consider the case $\psi=0$. Now all invariants of degrees 4 or 14
(or 18) vanish, but the condition $\psi=0$ itself yields the three
equations $A=B=C=0$ where $\psi=Ax^2+Bxy+Cy^2$.
Earlier, the choice $\psi \sim x^2$ used up some of the freedom given
by the group, but here we are free to choose a zero for the form,
and assume $h=0$. Again consider the four cases, this time with
$ag-6bf+15ce-10d^2$ zero instead of nonzero.
We have (iii) $f=1$, $h=g=e=d=c=b=0$, only variables $a,f$ left.
And (ii) $g=1$, $h=f=0$, $b+de = 0$, $3c+5e^2 = 0$,
$a+15ce-10d^2 = 0$, only variables $d,e$ left.
And by assumption $h = 0$ we are not in case (i).
That settles the case $\psi=0$.

Now assume $\psi \ne 0$ and take $\psi \sim x^2$.
In case (iii) only variables $a,b$ are left, and we are done.
In case (ii) only variables $a,d,e$ are left.
In case (i) only variables $c,d,e,f$ are left.
An invariant of degree 14 is
$(f.(f,\psi_2)_5,\psi^5)_{10} \sim (f.(f,\psi_2)_5,x^{10})_{10} =
-2afh^2+2ag^2h+7beh^2-7bfgh-5cdh^2-22cegh+27cf^2h+25d^2gh-45defh+20e^3h$.
In case (ii) this vanishes.
In case (i) this becomes (up to a constant) $18e^3-32def+9cf^2-cd$.
Another invariant of degree 14 is
$((\psi_2,\psi_3)_1,\psi^4)_8 \sim ((\psi_2,\psi_3)_1,x^8)_8$.
In case (ii) this becomes $de(26e^3-35d^2-10a)$ and we are reduced
to three pieces, each with only two variables.
In case (i) this becomes (up to a constant)
$70e^3f^4-120def^5+27cf^6+36e^5f-60de^3f^2+6ce^2f^3+3cdf^4+6d^2e^3+
18ce^4-8d^3ef-54cde^2f+33cd^2f^2+3c^2ef^2+cd^3-3c^2de+2c^3f$.
Both polynomials found are irreducible and hence have no common factor,
and we are reduced to a 2-dimensional situation. This proves $[4,14] \ge 3$.

\medskip
Finally, let us show that $[4,18] \ge 3$. The subcase $\psi=0$ was
handled already, so we can assume that $\psi \ne 0$ and take $\psi \sim x^2$.
Again only cases (i) and (ii) need to be considered.
Above we already considered the invariant $((\psi_1,\psi_2)_1,\psi^7)_{14}$
of degree 18. In case (ii) this yields $d=0$, leaving only the two
variables $a,e$. In case (i) we find $ef^2+de-cf = 0$.
Another invariant of degree 18 is
$(f.((f,\psi_2)_5,\psi_2)_2,\psi^6)_{12}$. In case (i) this yields
$70e^3f^3-120def^4+27cf^5-54e^5+210de^3f-200d^2ef^2-15ce^2f^2+
30cdf^3+15cde^2-25cd^2f-c^3 = 0$.
Both polynomials found are irreducible and hence have no common factor,
and we are reduced to a 2-dimensional situation. This proves $[4,18] \ge 3$.
\qed

\subsection{$n=8$}
For the octavic there there are nine basic invariants $i_d$ $(2 \le d \le 10)$.
There is a hsop with degrees 2, 3, 4, 5, 6, 7. The Poincar\'e series is
\begin{eqnarray*}
P(t) & = & 1 + t^2 + t^3 + 2t^4 + 2t^5 + 4t^6 + 4t^7 + 7t^8 + 8t^9 + \\
     && 12t^{10} + 13t^{11} + 20t^{12} + 22t^{13} + 31t^{14} + \cdots ~=~ \\
     & = & (1+t^8+t^9+t^{10}+t^{18}) / \prod_{d=2}^7 (1-t^d). \\
\end{eqnarray*}
\vskip -0.5cm
Given a finite sequence $(d_i)$, the {\em numerator} of $P(t)$
corresponding to this sequence is by definition
$P(t) \prod (1-t^{d_i})$. If $(d_i)$ is a subsequence of the
sequence of degrees of a hsop, then the corresponding numerator
has nonnegative coefficients. This rules out, e.g., the following
sequences $(d_i)$.

\medskip
\begin{tabular}{llll}
2, 2    & 2, 4, 4 & 3, 5, 5 & 5, 5, 5     \\
3, 3    & 2, 5, 5 & 4, 4, 4 & 2, 3, 7, 7  \\
\end{tabular}

\medskip
What is wrong with these sequences is that there just aren't enough
invariants of these degrees. More interesting are the cases where
there are enough invariants, but they cannot be chosen algebraically
independent.

\begin{Proposition}
A sequence of six integers larger than $1$ is the sequence of degrees
of a hsop for the octavic if and only if

(i) (`divisibility') at least three of them are even,
at least two are divisible by $3$, at least one has a factor $4$,
at least one a factor $5$, at least one a factor $6$,
and at least one a factor $7$, and moreover

(ii) (`nonnegativity') none of the eight sequences in the
above table occur as a subsequence, and moreover

(iii) (`algebraic independence')
there are no four elements in any of $\{2,3,6\}$, $\{2,4,5\}$, $\{2,4,7\}$,
and no five elements in any of $\{2,3,4,5,11\}$, $\{2,3,4,6,11\}$,
$\{2,3,4,7\}$, $\{2,3,4,8\}$, $\{2,3,4,9\}$, $\{2,3,5,6\}$, $\{2,3,6,7,11\}$.
\end{Proposition}
\Proof
We have
\[
I_2 = \langle i_2 \rangle,~~
I_3 = \langle i_3 \rangle,~~
I_4 = \langle i_2^2, i_4 \rangle,~~
I_5 = \langle i_2i_3, i_5 \rangle,~~
I_6 = \langle i_2^3, i_2i_4, i_3^2, i_6 \rangle,
\]
\[
I_7 = \langle i_2^2i_3, i_2i_5, i_3i_4, i_7 \rangle,~~
I_8 = \langle i_2^4, i_2^2i_4, i_2i_3^2, i_2i_6, i_3i_5, i_4^2, i_8 \rangle ,
\]
\[
I_9 = \langle i_2^3i_3, i_2^2i_5, i_2i_3i_4, i_2i_7,
i_3^3, i_3i_6, i_4i_5, i_9 \rangle ,
\]
\[
I_{11} = \langle i_2^4i_3, i_2^3i_5, i_2^2i_3i_4, i_2^2i_7,
i_2i_3^3, i_2i_3i_6, i_2i_4i_5, i_2i_9, i_3^2i_5, i_3i_4^2,
i_3i_8, i_4i_7, i_5i_6 \rangle .
\]

\medskip
We see that $V(\cup_{a\in A} I_a) = V(\{i_b \mid b \in B\})$
for $A$ and $B$ as in the table below.

\medskip
\begin{tabular}{cc|cc|cc}
$A$ & $B$ & $A$ & $B$ & $A$ & $B$ \\
\hline
2,3,6      & 2,3,6    & 2,3,4,6,11 & 2,3,4,6  & 2,3,5,6    & 2,3,5,6 \\
2,4,5      & 2,4,5    & 2,3,4,7 & 2,3,4,7     & 2,3,6,7,11 & 2,3,6,7 \\
2,4,7      & 2,4,7    & 2,3,4,8 & 2,3,4,8     & \\
2,3,4,5,11 & 2,3,4,5  & 2,3,4,9 & 2,3,4,9 & \\
\end{tabular}

\medskip
This shows that the given conditions are necessary.
For sufficiency, use induction. The basis of the induction
is provided by the 13 hsops constructed in the next proposition.
Given a sequence of six numbers satisfying the conditions, order
the numbers in such a way that the last is divisible by 7
and at least one of the last two is divisible by 5.
All restrictions concern numbers at most 11, so if we
split a number from the sequence into two parts each at least 12,
such that the divisibility conditions remain true for the two
resulting sequences, then by Lemma \ref{mainlemma} and induction
there exists a hsop with the given sequence as degree sequence.
This means that one can reduce the first four numbers modulo 12,
the fifth modulo 60, and the last modulo 420.
It remains to check a $24 \times 24 \times 24 \times 24 \times 72 \times 432$
box, and this is done by a small computer program. \qed

\begin{Proposition}
There are precisely $13$ minimal degree sequences of hsops in case $n = 8$,
namely

\[
\begin{array}{llll}
2,3,4,5,6,7   & 2,3,4,6,9,35  & 2,3,5,6,10,28 \\
2,3,4,5,8,42  & 2,3,4,7,8,30  & 2,3,5,9,12,14 \\
2,3,4,5,9,42  & 2,3,4,7,9,30  & 2,4,5,6,8,21  \\
2,3,4,5,10,42 & 2,3,4,8,9,210 & \\
2,3,4,6,8,35  & 2,3,5,6,9,28  & \\
\end{array}
\]
\end{Proposition}
\Proof
Minimality is immediately clear, so we only have to show existence.
Apply Dixmier's criterion.
As before we  have to show that for all $p$ and each subsequence
$d_1,...,d_p$ of one of these 13 sequences the inequality
$[d_1,...,d_p] \ge p$ holds.

We can save some work by observing that Shioda \cite{Sh} already
showed the existence of a hsop with degree sequence 2, 3, 4, 5, 6, 7.
It follows that $[d_1,...,d_p] \ge p$ when (at least) $p$ of the numbers
2, 3, 4, 5, 6, 7 divide some of the $d_i$.

For $p=1$, nothing remains to check.
% we need $[d] \ge 1$ for certain $d \ge 2$,
% and this follows from Proposition \ref{zeroh}.

For $p=2$, there only remains to show $[9] \ge 2$, and this follows
since there are two invariants of degree 9 without common factor,
for example $i_3i_6$ and $i_4i_5$.
% This settles $p=2$.

For $p=3$, we have to show $[8] \ge 3$, $[2,9] \ge 3$, $[5,9] \ge 3$,
$[7,9] \ge 3$, $[10] \ge 3$.

For $p=4$, we have to show $[3,8] \ge 4$, $[5,8] \ge 4$,
$[7,8] \ge 4$, $[4,9] \ge 4$, $[2,5,9] \ge 4$, $[6,9] \ge 4$,
$[2,7,9] \ge 4$, $[8,9] \ge 4$, $[3,10] \ge 4$, $[4,10] \ge 4$,
$[9,14] \ge 4$.

For $p=5$, we have to show
$[3,5,8] \ge 5$, $[6,8] \ge 5$, $[3,7,8] \ge 5$, $[4,5,9] \ge 5$,
$[4,6,9] \ge 5$, $[5,6,9] \ge 5$, $[4,7,9] \ge 5$, $[8,9] \ge 5$,
$[3,4,10] \ge 5$, $[6,10] \ge 5$, $[5,9,14] \ge 5$.

There are no conditions left to check for $p=6$.

Remain 27 conditions to check. Let $V[d_1,...,d_p]$ denote
the variety defined by all invariants of degrees $d_i$.
Split $V[9]$ into two parts depending on whether $i_2$ vanishes or not.
Where it does not vanish, all invariants of degrees 3, 5, 7 must vanish.
Hence $[5,9], [7,9] \ge [9] \ge \min([2,9], [3,5,7,9])$.
Split $[2,9]$ into two parts depending on whether $i_4$ vanishes or not.
The first part has $[2,3,4,9] \ge 3$, the second $[2,3,5,9] \ge 3$.
Hence $[9] \ge 3$.
Similarly, $[8] = [2,4,8] \ge \min([2,3,4,8],[2,4,5,8]) \ge 3$.
Finally, $[10] = [2,5,10] \ge \min([2,3,5,10],[2,3,7,10]) \ge 3$.
This settles $p=3$.

The same argument shows that
$[7,8], [2,7,9], [6,9], [3,10], [4,10], [9,14] \ge 4$
and $[5,9,14] \ge 5$.
% $[7,8] \ge \min([2,3,4,7],[2,4,5,7]) \ge 4$
% and $[2,7,9] \ge \min([2,3,4,7],[2,3,5,7]) \ge 4$
% and $[6,9] \ge \min([2,3,4,6],[2,3,5,6]) \ge 4$
% and $[3,10] \ge \min([2,3,4,5],[2,3,5,6]) \ge 4$
% and $[4,10] \ge \min([2,3,4,5],[2,4,5,7]) \ge 4$
% and $[9,14] \ge \min([2,3,4,7],[2,3,5,7]) \ge 4$.
% $[5,9,14] \ge \min([2,3,4,5,7],[2,3,5,6,7]) \ge 5$.

Since adding a single condition diminishes the dimension by at most one,
$[3,8] \ge 4$ follows from $[3,5,8] \ge 5$. (Given that $i_2$ vanishes
since $i_2^4$ has degree 8, the condition that all invariants of degree
5 vanish is equivalent to the requirement that $i_5$ vanishes.)
Similarly $[5,8] \ge 4$ and $[4,9] \ge 4$ and $[2,5,9] \ge 4$
follow from $[3,5,8] \ge 5$ and $[4,5,9] \ge 5$.
Trivially, $[8,9] \ge 4$ follows from $[8,9] \ge 5$.
This settles $p=4$, assuming the inequalities for $p=5$.

\bigskip
Remain 10 conditions to check:
$[3,5,8] \ge 5$, $[6,8] \ge 5$, $[3,7,8] \ge 5$, $[4,5,9] \ge 5$,
$[4,6,9] \ge 5$, $[5,6,9] \ge 5$, $[4,7,9] \ge 5$, $[8,9] \ge 5$,
$[3,4,10] \ge 5$, $[6,10] \ge 5$.

Equivalently, for each of the sets $A$, where $A$ is one of

\[\begin{array}{@{}lllll@{}}
\{2,3,4,5,8\},&
\{2,3,4,6,8\},&
\{2,3,4,7,8\},&
\{2,3,4,5,9\},&
\{2,3,4,6,9\},\\
\{2,3,5,6,9\},&
\{2,3,4,7,9\},&
\{2,3,4,8,9\},&
\{2,3,4,5,10\},&
\{2,3,5,6,10\},
\end{array}\]

\noindent
we must have $\dim V(\{i_a\mid a \in A\}) = 1$.

For example, we want $\dim V(i_2,i_3,i_4,i_5,i_8) = 1$.
Now $i_2,i_3,i_4,i_5$ form part of a hsop, so $V(i_2,i_3,i_4,i_5)$
is irreducible and has dimension 2. Moreover $i_8$ does not vanish
identically on $V(i_2,i_3,i_4,i_5)$ as we shall see, and it follows that
$\dim V(i_2,i_3,i_4,i_5,i_8) = 1$.

This argument works in all cases except that of $V(i_2,i_3,i_4,i_8,i_9)$
and shows that each of the claimed sequences of degrees with the
possible exception of 2, 3, 4, 8, 9, 210, is that of a hsop.
In particular, e.g. 2, 3, 4, 5, 8, 42 is the sequences of degrees of a hsop.
But now this argument also applies to $V(i_2,i_3,i_4,i_8,i_9)$:
$V(i_2,i_3,i_4,i_8)$ is irreducible of dimension 2 and $i_9$ does not
vanish identically on it, and it follows that $V(i_2,i_3,i_4,i_8,i_9)$
has dimension 1.

\medskip
It remains to check the ten conditions that say that $i_8$ does not vanish
on any of $V(i_2,i_3,i_4,i_5)$, $V(i_2,i_3,i_4,i_6)$, $V(i_2,i_3,i_4,i_7)$,
that $i_9$ does not vanish on any of $V(i_2,i_3,i_4,i_5)$,
$V(i_2,i_3,i_4,i_6)$, $V(i_2,i_3,i_5,i_6)$, $V(i_2,i_3,i_4,i_7)$,
$V(i_2,i_3,i_4,i_8)$, and that $i_{10}$ does not vanish on
$V(i_2,i_3,i_4,i_5)$ or $V(i_2,i_3,i_5,i_6)$.
Using Singular we computed the radical of the ideals $(i_2,i_3,i_4,i_5)$,
$(i_2,i_3,i_4,i_6)$, $(i_2,i_3,i_4,i_7)$, $(i_2,i_3,i_5,i_6)$ and
$(i_2,i_3,i_4,i_8)$ and checked the required facts.

(This shows that $i_8$, $i_9$ and $i_{10}$ do not vanish on the 2-dimensional
pieces mentioned. Note that these invariants do vanish on various
1-dimensional pieces. For example,
$i_8^2 \in (i_2,i_3,i_4,i_6,i_7)$, so that $i_8$ vanishes on
$V(i_2,i_3,i_4,i_6,i_7)$, and $i_8^5 \in (i_2,i_3,i_4,i_5,i_6)$, and
$i_{10}^2 \in (i_2,i_3,i_4,i_5,i_6)$ and
$i_9^3 \in (i_2,i_3,i_4,i_5,i_6) \cap (i_2,i_3,i_4,i_6,i_7) \cap
(i_2,i_3,i_5,i_6,i_7)$.)

% Use the explicit formulas given by Shioda \cite{Sh}.
% So far $i_d$ was an unspecified basic invariant of degree $d$. Let us now
% follow Shioda and pick definite invariants and call them $j_d$.
% 
% Define the covariants
% \[\begin{array}{llll}
% H = (f,f)_2, &
% g = (f,f)_4, &
% k = (f,f)_6, &
% h = (k,k)_2, \\
% m = (f,k)_4, &
% n = (f,h)_4, &
% p = (g,k)_4, &
% q = (g,h)_4
% \end{array}\]
% and invariants
% \[\begin{array}{lllll@{}}
% j_2 = (f,f)_8, &
% j_3 = (f,g)_8, &
% j_4 = (k,k)_4, &
% j_5 = (m,k)_4, &
% j_6 = (k,h)_4, \\
% j_7 = (m,h)_4, &
% j_8 = (p,h)_4, &
% j_9 = (n,h)_4, &
% j_{10} = (q,h)_4.
% \end{array}\]
%
% Modulo $(j_2,j_3,j_4,j_5)$ we have
% $j_8^2 \in (2^6/3^2 5)j_6j_{10} + (-2/5)j_7j_9$.
% Hence $j_8 \in \Rad (j_2,j_3,j_4,j_5,j_6)$.

\qed

\end{document}